\documentclass[12pt]{amsart}
\usepackage{amsfonts}
\begin{document}

\newtheorem{result}{Result}[section]
\newtheorem{thm}[result]{Theorem}
\newtheorem{prop}[result]{Proposition}
\newtheorem{lemma}[result]{Lemma}
\newtheorem{cor}[result]{Corollary}
\newtheorem{remark}[result]{Remark}

\newtheorem{defin}{Definition}
\newtheorem{example}{Example}
\newtheorem{problem}{Problem}

\def\dist{\mathop{\rm dist}}
\def\N{{\mathbb N}}
\let\<\langle
\let\>\rangle
\let\eps\varepsilon

\def\ind{\mathop{\rm ind}}
\def\span{\mathop{\rm span}}

\def\cal{\fam\symsymbols}
\def\Lc{{\cal L}}
\def\cop{{\cal C}}
\def\dop{{\cal D}}
\def\K  {{\cal K}}
\def\SS {{\cal SS}}
\def\SC {{\cal SC}}
\def\In {{\cal I}n}

\author[M.\ Gonz\'alez]{Manuel Gonz\'alez}
\address{Departamento de Matem\'aticas, Facultad de Ciencias,
Universidad de Cantabria, Avda. de los Castros s/n, 39071-Santander, Spain}
\email{manuel.gonzalez@unican.es}

\author[F.\ Le\'on-Saavedra]{Fernando Le\'on-Saavedra}
\address{Regional Mathematical Center of Southern Federal University, Rostov-on-Don, Russia \&
Department of Mathematics, University of C\'adiz, Avda. de la Universidad s/n,
11402-Jerez de la Frontera, Spain.}
\email{fernando.leon@uca.es}

\author[M.P.\ Romero de la Rosa]{Mar\'{\i}a del Pilar Romero de la Rosa}
\address{Department of Mathematics, University of C\'adiz, CASEM, Pol. R\'{\i}o San Pedro s/n,
11510-Puerto Real, Spain.}
\email{pilar.romero@uca.es}

\title[On $\ell_\infty$-Grothendieck subspaces]{On $\ell_\infty$-Grothendieck subspaces}

\thanks{The authors were supported by Ministerio de Ciencia, Innovaci\'on y Universidades (Spain), 
grants PGC2018-101514-B-I00,  PID2019-103961GB-C22, and by  Vicerrectorado de Investigaci\'on de 
la Universidad de C\'adiz.
This work was also co-financed by the 2014-2020 ERDF Operational Programme, and by the Department 
of Economy, Knowledge, Business and University of the Regional Government of Andalusia. 
Project reference: FEDER-UCA18-108415.\\
2010 Mathematics Subject Classification. Primary:  46A35; 46B20; 40H05.\\
Keywords: Grothendieck Banach space; $\ell_\infty$-Grothendieck subspace.}

\begin{abstract}
A closed subspace $S$ of $\ell_\infty$ is said to be a \emph{$\ell_\infty$-Grothendieck subspace} 
if $c_0\subset S$ (hence $\ell_\infty\subset S^{**}$) and every $\sigma(S^*,S)$-convergent sequence 
in $S^*$ is $\sigma(S^*,\ell_\infty)$-convergent.
Here we give examples of closed subspaces of $\ell_\infty$ containing $c_0$ which are or fail to 
be  $\ell_\infty$-Grothendieck.
\end{abstract}

\maketitle

\thispagestyle{empty}

\section{Introduction}

The $\ell_\infty$-Grothendieck subspaces (defined in the abstract; see also Definition \ref{ell_infty-Groth}) 
naturally emerge when some versions of Schur's Lemma for bounded multiplier convergent series are sharpened 
(see, e.g., \cite{Sw,AGP:07,AAGP:11,LRS-mdpi:20}).

Apart from $\ell_\infty$, only one example of $\ell_\infty$-Grothendieck subspace is given in the literature 
(see \cite[Remark 4.2]{AAGP:11}), using a result of \cite{Haydon:11}. 
This example is isomorphic to a $C(K)$ space with the Grothendieck property, and contains no subspaces 
isomorphic to $\ell_\infty$.

Here we prove that if $X$ is a Grothendieck Banach space and $M$ is a closed subspace of $X$ with $X/M$ 
separable then $M$ is a Grothendieck space. As a consequence, we derive that a closed subspace $S$ 
of $\ell_\infty$ containing $c_0$ is $\ell_\infty$-Grothendieck when the quotient $\ell_\infty/S$ is 
separable, and using the fact that $L_q(0,1)$ is isomorphic to a quotient of $\ell_\infty$ for 
$2\leq q<\infty$, we prove the existence of an uncountable family of pairwise non-isomorphic 
$\ell_\infty$-Grothendieck subspaces.
We also show that for each closed subspace $Y$ of $\ell_\infty$ which is a Grothendieck space and 
contains a subspace isomorphic to $c_0$, there exists a $\ell_\infty$-Grothendieck subspace 
isomorphic to $Y$.

On the other hand, we show that a closed subspace $S$ of $\ell_\infty$ containing $c_0$ is not 
$\ell_\infty$-Grothendieck when it is separable or, more generally, when the unit ball of $S^*$ 
is weak$^*$-sequentially compact.

\section{Preliminaries}

Let $X$ be a Banach space and let $M$ be a subspace of $X^{**}$ containing $X$. 
We say that a sequence $(x^*_n)$ in $X^*$ is $\sigma(X^*,M)$-convergent to $x^*$ if 
$(\langle x^{**},x^*_n\rangle)$ converges to $\langle x^{**},x^*\rangle$ for every $x^{**}\in M$.

A Banach space $X$ \emph{has weak$^*$ sequentially compact dual ball} (has W$^*$SC dual ball, for 
short) if every sequence in the unit ball of $X^*$ has a $\sigma(X^*,X)$-convergent subsequence. 
We refer to \cite[Chapter XIII]{Diestel:84} for information about this property. 
The next result gives some examples of spaces of this kind.

\begin{prop}\label{W*SCDB-sp}
A Banach space $X$ has W$^*$SC dual ball in the following cases:
\begin{enumerate}
  \item $X$ is separable;
  \item $X^*$ contains no copies of $\ell_1$;
  \item $X$ is isomorphic to the dual of a separable space containing no copies of $\ell_1$.
\end{enumerate}
\end{prop}
\begin{proof}
(1) is well-known \cite{Diestel:84}; (2) follows from Rosenthal's characterization of Banach spaces containing no copies of $\ell_1$ and the fact that each weakly Cauchy sequence in $X^*$ is weak$^*$-convergent; and (3) is a consequence of \cite[Theorem XIII.10]{Diestel:84}.
\end{proof}

A Banach space $X$ is \emph{Grothendieck} if every $\sigma(X^*,X)$-convergent sequence in $X^*$ is $\sigma(X^*,X^{**})$-convergent.

Obviously reflexive Banach spaces are Grothendieck. Moreover, it is not difficult to show that Grothendieck spaces with W$^*$SC dual ball are reflexive, it was proved in \cite{Grothendieck:53} that $\ell_\infty$ is a non-reflexive Grothendieck space (see \cite[Theorem VII.15]{Diestel:84}), and the class of Grothendieck spaces satisfies  the three-space property: If $M$ is a closed subspace of a Banach space $X$ and both $M$ and $X/M$ are Grothendieck, then so is $X$ (see  \cite[Corollary 2.6]{GO:86}).

The following result collects  some classical characterizations of Grothendieck spaces taken from
\cite[Chapter 5, Corollary 5]{Diestel:75}.

\begin{prop}\label{Groth-sp}
For a Banach space $X$, the following assertions are equivalent:
\begin{enumerate}
  \item $X$ is Grothendieck;
  \item every operator $T:X\to c_0$ is weakly compact;
  \item for each separable Banach space $Y$, every operator $T:X\to Y$ is weakly compact.
\end{enumerate}
\end{prop}

It easily follows from Proposition \ref{Groth-sp} that quotients of a Grothendieck space are also Grothendieck.

\section{Main results}

The following notion extends the classical one of Grothendieck space.

\begin{defin}\label{M-Groth}
Let $X$ be a Banach space and let $M$ be a vector subspace of $X^{**}$ containing $X$.
We say that $X$ is a \emph{$M$-Grothendieck space} if every $\sigma(X^*,X)$-convergent sequence
in $X^*$ is $\sigma(X^*,M)$-convergent.
\end{defin}

Obviously, the Grothendieck spaces are the $X^{**}$-Grothendieck spaces. Moreover, since $\sigma(X^*,X)$-convergent sequences are bounded, the $M$-Grothendieck spaces coincide with the $\overline{M}$-Grothendieck spaces, where $\overline{M}$ is the closure of $M$. So we could always assume in Definition \ref{M-Groth} that $M$ is a closed subspace.

We are interested in a concrete case of Definition \ref{M-Groth}. Let $S$ be a closed subspace of $\ell_\infty$ containing $c_0$, and let $j:c_0\to S$ be the inclusion map. Then we can identify $\ell_\infty$ with a subspace $j^{**}(c_0^{**})$ of $S^{**}$ containing $S$.

\begin{defin}\label{ell_infty-Groth}
Let $S$ be a closed subspace of $\ell_\infty$. We say that $S$ is a \emph{$\ell_\infty$-Grothendieck subspace} if it contains $c_0$ and each $\sigma(S^*,S)$-convergent sequence in $S^*$ is  $\sigma(S^*,\ell_\infty)$-convergent.
\end{defin}

Clearly, if $S$ is a closed subspace of $\ell_\infty$ that contains $c_0$ and $S$ is a Grothendieck space, then $S$ is a $\ell_\infty$-Grothendieck subspace.

The following result may be interesting on its own.

\begin{prop}
Let $X$ be a Grothendieck Banach space. If $M$ is a closed subspace of $X$ and $X/M$ is separable, then $M$ is a Grothendieck space.
\end{prop}
\begin{proof}
Let $S:M\to c_0$ be an operator. Since the space $c_0$ is separably injective \cite[Theorem 2.3]{ACCGM:16} and the quotient $X/M$ is separable, the operator $S$ admits an extension $T:X\to c_0$ \cite[Proposition 2.5]{ACCGM:16}, which is weakly compact by Proposition \ref{Groth-sp}. Then $S$ is weakly compact, and applying again Proposition \ref{Groth-sp} we conclude that $M$ is Grothendieck.
\end{proof}

As a consequence of the previous result, we obtain that ``big'' subspaces are $\ell_\infty$-Grothendieck subspaces.

\begin{cor}
Let $S$ be a closed subspace of $\ell_\infty$ containing $c_0$ such that $\ell_\infty/S$ is separable.
Then $S$ is a $\ell_\infty$-Grothendieck subspace.
\end{cor}

Let us see that there exists an uncountable family of pairwise non-isomorphic $\ell_\infty$-Grothendieck subspaces.

\begin{thm}\label{Lp-quot}
Let $2\leq p<\infty$.
\begin{enumerate}
\item There exists a closed subspace $N_p$ of $\ell_\infty$ containing $c_0$ such that the quotient $\ell_\infty/N_p$ is isomorphic to $L_p(0,1)$. Hence $N_p$ is a $\ell_\infty$-Grothendieck subspace.
\item If $2\leq r <\infty$, $p\neq r$, then the subspaces $N_p$ and $N_r$ are not isomorphic.
\end{enumerate}
\end{thm}
\begin{proof}
(1) Recall that $\ell_\infty$ is isomorphic to $L_\infty(0,1)$, which is the dual of $L_1(0,1)$. Let $q$ such that $1/p\, +\, 1/q=1$, hence $1<q\leq 2$.

By \cite[Corollary 2.f.5]{LT2}, there exists a closed subspace $M_q$ of $L_1(0,1)$ which is  isometrically isomorphic to $L_q(0,1)$. Therefore, by duality,
$$
M_q^*\equiv L_\infty(0,1)/M_q^\perp \equiv L_q(0,1)^*\equiv L_p(0,1).
$$

Let $U:L_\infty(0,1)\to\ell_\infty$ be a bijective isomorphism. By taking $N_p=U(M_q^\perp)$ we guarantee that $\ell_\infty/N_p$ is isomorphic to $L_p(0,1)$.

It remains to show that we can choose $N_p$ containing $c_0$. This is a consequence of the fact that $\ell_\infty/c_0$ has a quotient isomorphic to $\ell_\infty$. So we can take as $N_p$ the kernel of a composition of surjective operators like the following one:
$$
\ell_\infty\to\ell_\infty/c_0\to \ell_\infty\to L_p(0,1).
$$

(2) Let $2\leq p,r<\infty$, and assume that there exists a bijective isomorphism $T:N_p\to N_r$. Since both $\ell_\infty/N_p$ and $\ell_\infty/N_r$ are reflexive, by \cite[Theorem 2.f.12]{LT1} there exists an extension $\hat T :\ell_\infty\to\ell_\infty$ of $T$ which is a Fredholm operator; i.e. the range $R(\hat T)$ is closed and both the kernel $N(\hat T)$ and $\ell_\infty/R(\hat T)$ are finite dimensional. Then $\hat T$ induces a Fredholm operator $S:\ell_\infty/N_p\to \ell_\infty/N_r$, implying that $L_p(0,1)$ and $L_r(0,1)$ are isomorphic. Hence $p=r$, and the proof is done.
\end{proof}

Every quotient of $L_p(0,1)$ ($2\leq p<\infty$) is also a quotient of $\ell_\infty$, and we can assume as before that the kernel of the quotient map contains $c_0$, so it provides another example of $\ell_\infty$-Grothendieck subspace. In particular, we could have formulated Theorem
\ref{Lp-quot} with $\ell_p$ instead of $L_p(0,1)$.

The next result provides additional examples.

\begin{prop}
Let $Y$ be a closed subspace of $\ell_\infty$ which is a Grothendieck space and contains a subspace isomorphic to $c_0$. Then there exists a $\ell_\infty$-Grothendieck subspace isomorphic to $Y$.
\end{prop}
\begin{proof}
Let $M$ be a closed subspace of $Y$ isomorphic to $c_0$ and let $T:M\to c_0$ be a bijective isomorphism. Since both $\ell_\infty/M$ and $\ell_\infty/c_0$ are non-reflexive, by \cite[Theorem 2.f.12]{LT1} there exists an extension $\hat T : \ell_\infty \to\ell_\infty$ of $T$ which is a bijective isomorphism. Hence $\hat T(Y)$ is a $\ell_\infty$-Grothendieck subspace isomorphic to $Y$.
\end{proof}

A remarkable example of Grothendieck space obtained by Bourgain \cite{Bourgain:83} is the space $H^\infty$ of bounded analytic functions on the unit disc, which is not isomorphic to a $C(K)$ space, not even isomorphic to a $\mathcal{L}_\infty$-space. Moreover, it was proved in \cite[Corollary 10]{GonzalezG:13} that the projective tensor product $\ell_\infty \widehat\otimes_\pi\ell_p$ is Grothendieck for $2<p<\infty$.

Since both spaces $H^\infty$ and $\ell_\infty\widehat\otimes_\pi \ell_p$ contain a subspace isomorphic to $c_0$ and they are isomorphic to dual spaces of separable spaces $(L_1/H_0^1)^*$ and $(\ell_1\widehat \otimes_\varepsilon \ell_p^*)^*$, hence they embed in $\ell_\infty$, we get the following fact.

\begin{cor}
There exist $\ell_\infty$-Grothendieck subspaces which are isomorphic to $H^\infty$ and  $\ell_\infty\widehat \otimes_\pi \ell_p$ for $2<p<\infty$.
\end{cor}

All known examples of $\ell_\infty$-Grothendieck subspace are Grothendieck spaces. So the following question arises:

\begin{problem}
Is it possible to find an example of $\ell_\infty$-Grothendieck subspace which is not a Grothendieck space?
\end{problem}

To study this problem we would need a good characterization of $\ell_\infty$-Grothendieck subspaces, which we do not have yet.

Next we show that ``small'' subspaces are not $\ell_\infty$-Grothendieck subspaces.

\begin{prop}\label{prop:noLG}
Let $S$ be a closed subspace of $\ell_\infty$ containing $c_0$.
If $S$ has W$^*$SC dual ball, then $S$ is not a $\ell_\infty$-Grothendieck subspace.
\end{prop}
\begin{proof}
Let $j:c_0\to S$ be the inclusion. Then $j^*:S^*\to c_0^*$ is surjective, and we can select a bounded sequence $(x^*_n)$ in $S^*$ such that $j^*x^*_n = e^*_n$ for each $n\in\N$, where $(e^*_n)$ is the unit vector basis of  $\ell_1\equiv c_0^*$.

Since $S$ has W$^*$SC dual ball, $(x^*_n)$ has a $\sigma(S^*,S)$-convergent subsequence. Thus the proof is finished if we show that $(x^*_n)$ has no $\sigma(S^*,\ell_\infty)$-convergent subsequence.

Indeed, let $(x^*_{n_k})$ be a subsequence, and  recall that $j^{**}:c_0^{**}\equiv\ell_\infty\to S^{**}$ is the inclusion.
We take $z=(a_i)\in\ell_\infty$ with $a_i = 1$ for $i=n_{2k}$ ($k\in\N$) and $a_i=-1$ otherwise. Then
$$
\langle j^{**}z,x^*_{n_k}\rangle = \langle z, j^*x^*_{n_k}\rangle =\langle z,e^*_{n_k}\rangle =(-1)^k,
$$
hence  $(x^*_{n_k})$ is not  $\sigma(S^*,\ell_\infty)$-convergent.
\end{proof}

Proposition \ref{prop:noLG} applies in the following cases:
\smallskip

(1) $S$ is a separable closed subspace of $\ell_\infty$ containing $c_0$.

(2) $S=\overline{c_0+M}$, where $M$ is a non-separable subspace of $\ell_\infty$ containing no copies of $\ell_1$ and isomorphic to a dual separable space. The space $S$ has W$^*$SC dual ball because $M$ has W$^*$SC dual ball (Proposition \ref{W*SCDB-sp}) and there is an injective operator with dense range $T:c_0\times M\to S$; see \cite[Chapter XIII]{Diestel:84}.
For example, we can take $M$ isomorphic to the dual of the James tree space $JT$ \cite{James}.
\medskip

\noindent
{\bf Acknowledgements.}
We  thank our colleagues of the Department of Mathematics, University of Cantabria for their support during a research stay.

\end{document}